\documentclass[12pt]{article}
\usepackage{amssymb,a4}
\usepackage{amsthm}

\def\N{\mathbb{N}}
\def\z{\mathbb{Z}}
\def\c{\mathbb{C}}
\def\dim{\hbox{dim}}

\def\a{\alpha}
\def\b{\beta}
\def\F{\mathbb C}

\def\ll{\lambda}

\def\mod{\hbox{mod}}

\newfont{\df}{eufm10}

\def\ll{\lambda}

\def\ot{\otimes}

\def\de{\delta}
\def\dim{\hbox{\rm dim}\,}

\def\Vir{\hbox{\rm Vir}}

\def\ot{\otimes}

\hoffset \voffset \oddsidemargin=60pt \evensidemargin=45pt
\title{\bf  Classification of Harish-Chandra
modules over the $W$-algebra $W(2,2)$
\thanks{This work is supported in part by NSF of China (No.
10671027, 10701019 and 10571119).}}
\author{Dong Liu\\ Department of Mathematics, Huzhou Teachers College\\ Zhejiang Huzhou, 313000, China\\
Linsheng Zhu\\ Department of Mathematics, Changshu Institute of
Technology\\ Jiangsu Changshu, 215500, China}
\date{ }

\begin{document}
\maketitle

\begin{abstract}  In this paper, we classify all irreducible weight modules with finite
dimensional weight spaces over the $W$-algebra $W(2, 2)$.
Meanwhile, all indecomposable modules with one dimensional weight
spaces over the $W$-algebra $W(2, 2)$ are also determined.
\end{abstract}



\smallskip\bigskip

\newtheorem{theo}{Theorem}[section]
\newtheorem{defi}[theo]{Definition}
\newtheorem{lemm}[theo]{Lemma}
\newtheorem{coro}[theo]{Corollary}
\newtheorem{prop}[theo]{Proposition}

\section{Introduction}

The $W$-algebra $W(2, 2)$ was introduced in \cite{DZ} for the
study of classification of vertex operator algebras generated by
weight 2 vectors. By definition, {\it the $W$-algebra $W(2, 2)$}
is the Lie algebra ${\cal L}$ with  $\c$-basis $\{L_m,I_m,C,
C_1|\, m \in \mathbb{Z}\}$ subject to the following relations.

\begin{defi} The {\bf $W$-algebra} ${{\cal L}}=W(2, 2)$ is
a Lie algebra over $\c$ (the  field of complex numbers) with the
basis
$$\{x_n,I(n),C, C_1 | n \in \z\}$$
and the Lie bracket given by
$$[x_n,x_m]=(m-n)x_{n+m}+\delta_{n,-m}\frac{n^3-n}{12}C,\eqno (1.1)$$
$$[x_n,I(m)]=(m-n)I(n+m)+\delta_{n,-m}\frac{n^3-n}{12}C_1, \eqno (1.2)$$
$$[I(n),I(m)]=0, \eqno (1.3)$$
$$[{\cal L},C]=[{\cal L},C_1]=0. \eqno (1.4)$$
\end{defi}

The $W$-algebra $W(2,2)$ can be realized from the semi-product of
the Virasoro algebra Vir and the Vir-module ${\cal A}_{0, -1}$ of
the intermediate series in \cite{OR}. In fact, let $W=\c\{x_m\mid
m\in\z\}$ be the Witt algebra (non-central Virasoro algebra) and
$V=\c\{I(m)\mid n\in\z\}$ be a $W$-module with the action
$x_m\cdot I(n)=(n-m)I(m+n)$, then $W(2, 2)$ is just the universal
central extension of the Lie algebra $W\ltimes V$ (see \cite{OR}
and \cite{GJJ}).  The $W$-algebra $W(2,2)$ studied in \cite{DZ} is
the restriction for $C_1=C$ of $W(2,2)$ in our paper.

The $W$-algebra $W(2,2)$ can be also realized from the so-called
{\it loop-Virasoro algebra} (see \cite{GLZ}). Let $\c[t, t^{-1}]$
be the Laurents polynomial ring over $\c$, then the loop-Virasoro
algebra $\tilde{VL}$ is the universal central extension of the
loop algebra $\Vir\ot\c[t^1, t^{-1}]$ and $W(2,
2)=\tilde{VL}/\c[t^2]$.

 The $W$-algebra
$W(2,2)$ is an extension of the Virasoro algebra and is similar to
the twisted Heisenberg-Virasoro algebra (see \cite{ADKP}).
However, unlike the case of the later, the action of $I(0)$ in
$W(2, 2)$ is not simisimple, so its representation theory is very
different from that of the twisted Heisenberg-Virasoro algebra in
a fundamental way.

The $W(2,2)$ and its highest weight modules enter the picture
naturally during our discussion on $L(1/2,0)\otimes L(1/2,0).$ The
$W$-algebra $W(2,2)$ is an extension of the Virasoro algebra and
also has a very good highest weight module theory (see Section 2).
Its highest weight modules produce a new class of vertex  operator
algebras. Contrast to the Virasoro algebra case, this class of
vertex operator algebras are always irrational.

The present paper is devoted to determining all irreducible weight
modules with finite dimensional weight spaces over  ${\cal L}$
from the motivations in \cite{LZ} and \cite{LJ}. More precisely we
prove that there are two different classes of them. One class is
formed by simple modules of intermediate series, whose weight
spaces are all $1$-dimensional; the other class consists of the
highest(or lowest) weight modules.

The paper is arranged as follows. In Section 2, we recall some
notations and collect known facts about irreducible,
indecomposable modules over the classical Virasoro algebra. In
Section 3, we determine  all irreducible (indecomposable) weight
modules of intermediate series over ${\cal L}$, i.e., irreducible
(indecomposable) weight modules with  all $1$-dimensional weight
spaces. In Section 4, we determine all irreducible uniformly
bounded weight modules over ${{\cal L}}$ which turn out to be
modules of intermediate series. In Section 5, we obtain the main
result of this paper: the classification of irreducible weight
${{\cal L}}$-modules with finite dimensional weight space. As we
mentioned, they are irreducible highest, lowest weight modules, or
irreducible modules of the intermediate series.

\section{Basics}

In this section, we collect some known facts for later use.

For any $e\in\c$, it is clear that
$$[x_n+neI(n),\; x_m+meI(m)]=(m-n)(x_{n+m}+(m+n)eI(n+m)),\,\, \forall n \neq -m,$$
$$[x_n+neI(n),\;
x_{-n}+(-n)eI(-n)]=-2nx_0+\frac{n^3-n}{12}C.$$ So $\{x_n+enI(n),
C| n\in \z\}$ spans a subalgebra $\Vir[e]$ which is isomorphic to
the classical Virasoro algebra. In many cases, we shall simply
write $\Vir[0]$ as $\Vir$.

\par
 Introduce a $\z$-grading on ${{\cal L}}$ by defining the degrees: deg $x_n$=deg $I(n)$=$n$ and deg
$C=0$. Set
 \vskip 5pt
$${\cal L}_+=\sum_{n\ge 0}(\c x_n+\c I(n)),\,\,{\cal L}_-=\sum_{n\le 0}(\c x_n+\c I(n)),$$
and
$${\cal L}_{0}=\c x_0+\c I(0)+\c C.$$

\vskip 5pt An ${\cal L}$-module $V$ is called a weight module if
$V$ is the sum of all its weight spaces. For a weight module $V$
we define
$$\hbox{Supp}(V):=\bigl\{\lambda\in \c \bigm|V_\lambda\neq
0\bigr\},$$ which is generally called the weight set (or the
support) of $V$.

\par
 A nontrivial weight ${\cal L}$-module V is called a {weight module of intermediate
 series} if V is indecomposable and any weight spaces of V is one dimensional.
\par
A  weight ${\cal L}$-module V is called a {highest} (resp.
{lowest) weight module} with {highest weight} (resp. {highest
weight}) $\lambda\in \c$, if there exists a nonzero weight vector
$v \in V_\lambda$ such that

 \vskip 5pt 1) $V$ is generated by $v$ as  ${\cal L}$-module;

2) ${\cal L}_+ v=0 $ (resp. ${\cal L}_- v=0 $).
 \vskip 5pt

\noindent{\bf Remark.} For a highest (lowest) vector $v$ we always
suppose that $I_0v=c_0v$ for some $c_0\in \c$ although the action
of $I_0$ is not semisimple.

Obviously, if $M$ is an irreducible weight ${\cal L}$-module, then
there exists $\lambda\in\c$ such that ${\rm
Supp}(M)\subset\lambda+\z$. So $M$ is a $\z$-graded module.

If, in addition, all weight spaces $M_\ll$ of a weight ${\cal
L}$-module $M$ are finite dimensional, the module is called a {\it
Harish-Chandra module}. Clearly a highest (lowest) weight module
is a Harish-Chandra module.

Let $U:=U(\cal L)$ be the universal enveloping algebra of  ${\cal
L}$. For any $\lambda,c$ $\in \c$, let $I(\lambda,c, c_0, c_1)$ be
the left ideal of $U$ generated by the elements $$
\bigl\{x_i,I(i)\bigm|i\in \N \bigr\}\bigcup\bigl\{x_0-\lambda
\cdot 1, C-c\cdot 1, I_0-c_0\cdot 1, C_1-c_1\cdot1\bigr\}. $$ Then
the Verma module with the highest weight $\lambda$ over ${\cal L}$
is defined as
$$M(\lambda, c, c_0, c_1):=U/I(\lambda, c, c_0, c_1).$$
It is clear that $M(\lambda, c, c_0, c_1)$ is a highest weight
module over ${\cal L}$ and   contains a unique maximal submodule.
Let $V(\lambda, c, c_0, c_1)$ be the unique irreducible quotient
of $M(\lambda, c, c_0, c_1)$.

The following result was given in \cite{DZ}.
\begin{theo} \cite{DZ} The Verma module $M(\lambda, c, c_0, c_1)$ is irreducible if and only if
$\frac{m^2-1}{12}c_1+2c_0\ne 0$ for any nonzero integer $m$.
\end{theo}

In the rest of this section, we recall some known facts about
weight representations of the classical Virasoro algebra which can
be considered as a subalgebra of $\cal L$:
$$\Vir:= \hbox{span} {\{x_n, C | n \in \z\}}.$$

Kaplansky-Santharoubane \cite{KS} in 1983 gave a classification of
$Vir$-modules of the intermediate series. There are three families
of {\it indecomposable modules of the intermediate series} (i.e
nontrivial indecomposable weight modules with each weight space is
at most one-dimensional) over the Virasoro algebra. They are
Vir-modules "without central charges".

(1) ${\mathcal A}_{a,\; b}=\sum_{i\in\z}\c v_i$:
$x_mv_i=(a+i+bm)v_{m+i}$;

(2) ${\mathcal A}(a)=\sum_{i\in\z}\c v_i$: $x_mv_i=(i+m)v_{m+i}$
if $i\ne 0$, $x_mv_0=m(m+a)v_{m}$;

(3) ${\mathcal B}(a)=\sum_{i\in\z}\c v_i$:  $x_mv_i=iv_{m+i}$ if
$i\ne -m$, $x_mv_{-m}=-m(m+a)v_0$,  for some $a, b\in\c$.

When $a\notin\z$ or $b\ne0, 1$, it is well-known that the module
${\cal A}_{a,\; b}$ is simple.  In the opposite case the module
contains two simple subquotients namely the trivial module and
$\c[t, t^{-1}]/\c$. Denote the nontrivial simple subquotients of
${\mathcal A}_{a,\; b}$, ${\mathcal A}(a)$, ${\mathcal B}(a)$ by
${\mathcal A}_{a, \; b}'$, ${\mathcal A}(a)'$, ${\mathcal B}(a)'$
respectively. They are all Harish-Chandra modules of the
intermediate series over the Virasoro algebra. (These facts appear
in many references, for example in \cite{SZ1}). We shall
 use $T$ to denote the $1$-dimensional trivial module,
 use $V'(0,0)$ to denote the unique proper nontrivial
submodule of $V(0,1)$ (which is irreducible).

An indecomposable module $V$ over $\Vir$ is said to be an {\it
extension of the $\Vir$-module $W_1$ by the $\Vir$-module $W_2$}
if $V$ has a submodule isomorphic to $W_1$ and $V/W_1\simeq W_2$.

\begin{theo}\cite{MP2}  Let $Z$ be an indecomposable weight $\Vir$-module
with weight spaces of dimension  less than or equal to two.  Then
$Z$ is  one of the following:

1) an intermediate series module;

2) an extension of $A_{\a,\b}$, $A(a)$ or $B(a)$ by themselves;

3) an extension of $A_{\a,\b_1}$ by $A_{\a,\b_2}$, where
$\b_1-\b_2=2,3,4,5,6$;

4) an extension of $A_{0, \b}$ by $W$, where $\b=2,3,4,5$ and $W$
is one of: $A_{0,0}$, $A_{0,1}$, $A_{0,0}'$, $A_{0,0}'\oplus T$,
$T$, $A(a)$ or $B(a)$;

5) an extension of  $T$ by   $A_{0,0}$ or $A(a)$, ;

6) an extension of $A_{0,0}'$ by  $A_{0,0}$ or $A(1)$;

7) an extension  $A_{0,1}$ by $A_{0,0}$ or $B(0)$;

8) an extension of $A_{0,0}$ by   $A_{0,1}$ or $A(a)$;

9) an extension of $A(0)$ by  $A_{0,0}$;

10) an extension of $B(a)$ by $A_{0,1}$;

11) the contragredient extensions of the previous ones;

\noindent where $\a,\b,\b_1,\b_2,a\in\c$.
\end{theo}

\noindent{\bf Remark.} In the above list, there are some
repetitions, and not all of them can occur.

\begin{theo}\cite{MP2} There are exactly two indecomposable
extensions $V=\hbox{span}\{v_{\a+i},\break v'_{\a+i}|i\in\z\}$ of
$A_{\a,0}$ by $A_{\a,0}$ ($\a \notin \z$) given by the actions $
x_iv_{\a+n}=(\a+n)v_{\a+n+i}, \forall i\in\z$ and

\noindent a) $x_iv_{\a+n}'=(\a+n)v_{\a+n+1}'-iv_{\a+n+1}$, for all
$i,n\in\z$; or

\noindent  b)  $x_1v_{\a+n}'=(\a+n)v_{\a+n+1}'$,
$x_{-1}v_{\a+n}'=(\a+n)v_{\a+n-1}'$,
 $x_2v_{\a+n}'=(\a+n)v_{\a+n+2}'+\frac{1}{(\a+n+2)(\a+n+1)}v_{\a+n+2}$,
$x_{-2}v_{\a+n}'=(\a+n)v_{\a+n-2}'-\frac{1}{(\a+n-2)(\a+n-1)}v_{\a+n-2},$
where $\{v_{\a+n},v'_{\a+n}\}$ forms a base of $V_{\a+n}$ for all
$n\in \z$.
\end{theo}

We also need the following result from \cite{MP1}.
\begin{theo}\cite{MP1}   Suppose that V is a
weight Vir-module with finite dimensional weight spaces. Let $M^+$
(resp. $M^-$ ) be the maximal submodule of V with upper (resp.
lower) bounded weights. If V and V* (the contragredient module of
V) do not contain trivial submodules, there exists a unique
bounded submodule $B$ such that $V=B\bigoplus M^+ \bigoplus M^-$.
\end{theo}

\section{Irreducible weight modules with weight multiplicity
one}

In this section we determine all irreducible and indecomposable
weight modules over ${\cal L}$ with weight multiplicity one.
\par
\vskip 2pt

Let $V=\oplus\c v_i $ be a $\z$-graded ${\cal L}$-module. Then the
action of $I(0)$ is semisimple and we can suppose that
$I(0)v_i=\ll_iv_i$ for some $\ll_i\in\c$. Moreover, by $[I(0),
x_m]=mI(m)$ we have $mI(m)v_i=(\ll_{m+i}-\ll_i)x_mv_i$. Hence $V$
is an irreducible (indecomposable) ${\cal L}$-module if and only
if $V$ is an irreducible (indecomposable) Vir-module.

Denoted ${\mathcal A}_{a,\; b,\; 0}$, ${\mathcal A}(a, 0)$,
${\mathcal B}(a, 0)$ by the $\cal L$-module induced from the
Vir-module ${\mathcal A}_{a,\; b}$, ${\mathcal A}(a)$, ${\mathcal
B}(a)$ with the trivial actions of $I(n)$ for any $n\in\z$,
respectively. Moreover we also denote the nontrivial simple
subquotients of ${\mathcal A}_{a,\; b}$, ${\mathcal A}(a)$,
${\mathcal B}(a)$ by ${\mathcal A}_{a, \; b}'$, ${\mathcal
A}'(a)$, ${\mathcal B}'(a)$ respectively. Clearly, ${\mathcal
A}'(a, 0)\cong {\mathcal B}'(a, 0)\cong {\mathcal A}_{0,\; 0,\;
0}'$
 Now we shall prove that the above three kinds
modules are all indecomposable modules with weight multiplicity
one.

\begin{lemm}\label{L1}
 Let $V=\sum_{i\in\z}\F v_i$ be an $\z$-graded ${\cal L}$-module
such that $x_nv_i=(a+i+bn)v_{i+n}$ for all $n, i\in\z$ and some
$a, b\in\c$. If $a+bn\ne0$ for any $n\in\z$, then $I(m)v_i=0$.

\end{lemm}
\noindent{\bf Proof.} Since $V$ is a module of the Virasoro
algebra Vir$=\oplus_{m\in\z} L(m)$, it is clear that $C=0$ (cf.
[SZ] for example). Suppose that $I(n)v_t=f(n, t)v_{n+t}$ for all
$n, t\in\z$. From $[L(n), I(m)]=(m-n) I(n+m)+{1\over 12}\de_{n+m,
0}(n^3-n)C_1$, we see that
$$f(m, t)L(n)v_{t+m}-f(m, n+t)(a+t+bn)v_{n+t+m}$$$$=
(m-n)f(n+m, t)v_{n+m+t}+{1\over 12}\de_{n+m, 0}(n^3-n)C_1v_{t}.
\eqno(3.1)$$ In this case,
$$f(m, t)(a+t+m+bn)-f(m, n+t)(a+t+bn)=(m-n) f(n+m, t)+{1\over 12}\de_{n+m, 0}(n^3-n)C_1. \eqno(3.2)$$

Let $t=0$ in (3.2), then
$$f(m, n)(a+bn)=
(a+m+bn)f(m, 0)-((m-n)f(n+m, 0)+{1\over 12}\de_{n+m,
0}(n^3-n)C_1).\eqno(3.3)$$

 Let $m=n$ in (3.2),
then
$$f(n, t)(a+(b+1)n+t)-f(n, n+t)(a+bn+t)=0. \eqno(3.4)$$

Let $t=0$ in (3.4), then
$$f(n, 0)(a+(b+1)n)-f(n, n)(a+bn)=0. \eqno(3.5)$$

Let $t=-n$ in (3.4), then
$$f(n, -n)(a+bn)=f(n, 0)(a+(b-1)n). \eqno(3.6)$$

Setting $n=-m$ in (3.3), then $$(a-m b)f(m,-m)=(a+m-bm)f(m,0)-2m
f(0,0)-{1\over 12}(m^{3}-m)C_1.\eqno(3.7)$$

From (3.6) and (3.7) we obtain that
$$
f(m, 0)={a+bm\over a}(f(0, 0)-{1\over 24}(m^2-1)C_1), \quad m\ne0.
\eqno(3.8)$$

Setting $t=n$ and $m=0$ in (3.2), we have
$$(f(0, 2n)-f(0, n))(a+(b+1)n)=nf(n, n). \eqno(3.9)$$

Let $m=0$ in (3.3), then
$$f(0, n)(a+bn)=
(a+bn)f(0, 0)+nf(n, 0).\eqno(3.10)$$

From (3.10) and (3.8) we have
$$
f(0,n)={a+n\over a}f(0, 0)-{1\over 24a}(n^3-n)C_1.\eqno(3.11)$$
Applying (3.11) to (3.9), we have
$$f(n,n)={a+(b+1)n\over a}(f(0,
0)-{1\over 24}(7n^2-1)C_1).\eqno(3.12)$$

Combining (3.5) and (3.11) we obtain $C_1=0.$

Therefore  $$ f(m, 0)={a+bm\over a}f(0, 0),$$ for all $m\in\z$. By
(3.3) we obtain that $$f(m, n)= {a+bm+n\over a}f(0,
0)=c(a+bm+n),\eqno(3.13)$$ for all $m, n\in\z$ and some
$c\in{\mathbb F}$.

However, by $[I(m), I(n)]=0$ we have $c=0$. So $f(m, n)=0$ for any
$m, n\in\z$. \hfill $\rule[-.23ex]{1.0ex}{2.0ex}$

\noindent{\bf Remark.} In the following  cases, we  can also
deduce that $C_1=C_1=0$  as in Lemma 7.3. So in the following
discussions, we always assume that $C_1=C_1=0$.

\begin{lemm}\label{L2} If $V\simeq {\cal A}_{\a,\;\b}$ as Vir-module,
 then $V\simeq {\cal A}_{\a,\;\b,\;0}$ as ${\cal L}$-module.
\end{lemm} \noindent{\bf Proof.}
 ({\bf I}.1) Suppose that
$a\not\in\z$ and $a+bn\ne 0$ for all $n\in\z$, then $f(m, n)=0$
for all $m, n\in\z$ by Lemma \ref{L1}.

({\bf I}.2) $a\not\in\z$ and $a=bp$ for some $p\in
\z\backslash\{0\}$. So $b\neq 0,1$. Therefore $f(m,n)=0$ if
$n+p\neq 0$ by (3.4).
 It follows from (3.4)
that $$f(0, -p)=0.\eqno(3.14)$$ Setting $i=-p$, $m=0$ in (3.2) and
using (3.14, we have
$$f(n,-p)=0, \quad n\ne0.$$
 Therefore $f(m, n)=0$ for all
$m,n\in\z$. \hfill $\rule[-.23ex]{1.0ex}{2.0ex}$

({\bf I}.3) $a\in\z$. Since ${\cal A}_{a, b}\cong {\cal A}_{0,
b}$, so we can suppose that $x_nv_i=(i+bn)v_{n+i}$ for all $n,
i\in\z$.

 ({\bf I}.3.1)\ $b\ne 0, 1$, then $a+bn\ne 0$ for all $n\ne 0$. So $f(m,
n)=0$ by Lemma \ref{L1}. Therefore $V$ is isomorphic to ${\cal
A}_{0,b,0}$.

({\bf I}.3.2) \ $b=1$. In this case (3.3) still holds and (3.2)
and (3.3) becomes
$$f(m, i)(i+m+n)-f(m, n+i)(i+n)=(m-n)f(n+m, i)\eqno(3.15)$$
and
$$
n f(m,n)=(n-m)f(n+m,0)+(n+m)f(m,0)\eqno(3.16)$$

Let $m=i=0$ in (3.15) we have
$$f(0, n)=f(n, 0)+f(0, 0), \quad n\ne0.\eqno(3.17)$$

Replacing $t$ by $n$ and letting $n=-m$ in (3.15), we have $$n
f(m,n)=(n-m)f(m,n-m)+2m f(0,n).\eqno(3.18)$$ Replacing $n$ by
$n-m$ in (3.16), we have
$$(n-m)f(m,n-m)=(n-2m)f(n,0)+nf(m,0).\eqno(3.19)$$
From (3.16)-(3.19), we have
$$
(n-m)f(n+m,0)=nf(n,0)-mf(m,0)+2mf(0,0).\eqno(3.20)$$ Let $n=0$ in
(3.20), we have $f(0,0)=0$. So (3.20) becomes
$$(n-m)f(n+m,0)=nf(n,0)-mf(m,0).$$

 We deduce
that there exists $c, d\in\F$ such that
$$
f(n,0)=c+dn,\eqno(3.21)$$ for all $n\in\z$ and $n\ne0$. Applying
(3.20) to (3.16) and using (3.21), we have
$$
f(m,n)=f(0,0), \ n\neq 0.$$ Therefore
  $$f(m, n)=2d+c(m+n), \quad n\ne0.$$

By (3.15) we have $d=0$.

Therefore $$f(m, n)=c(m+n).$$ By $[I(m), I(n)]=0$ we have $c=0$.
So $f(m, n)=0$ for any $m, n\in\z$.

 ({\bf I}.3.3) \ $b=0$. (3.2) becomes
$$
(i+m)f(m,i)-tf(m,n+i)=(m-n) f(m+n,i).\eqno(3.22)$$ Let $i=0$ and
$n=-m$ in (3.22), then
$$f(m,0)=f(0,0).$$
Let $i=0$ and $n=m$ in (3.22), then
$$f(m,0)=0, m\ne0.$$
Hence $$f(n,0)=0, \forall n\in\z.\eqno(3.23)$$

Setting $i=-n$ in (3.22) and using (3.23), we have
$$
f(m,-n)=f(m+n,-n),\quad m\ne n.\eqno(3.24)$$

Let $i=1 $ in (3.22) then
$$f(m, n+1)=(m+1)f(m, 1)-(m-n)f(m+n, 1).\eqno(3.25)$$
Setting $n=-1 $ in (3.25) and using (3.23), we have
$$f(m, 1)=f(m-1, 1), m\ne -1.\eqno(3.26)$$
So $f(m, 1)=c_1$ for some $c_1\in{\mathbb F}$ and for any $m\ge
-1$. $f(m, 1)=c_2$ for some $c_2\in{\mathbb F}$ and for any
$m\le-2$.

Hence (3.25) becomes
$$f(m, n)=c_1n,  \quad m\ge-1; \quad f(m, n)=c_2n,  \quad m\le-2.\eqno(3.27)$$
By (3.24) we have $c_1=c_2=c$ for some $c_1\in{\mathbb F}$.

Therefore $f(m, n)=nc$ for any $m, n\in\z$.
 By $[I(m), I(n)]=0$ we have $c=0$. So $f(m, n)=0$ for any
$m, n\in\z$. \hfill $\rule[-.23ex]{1.0ex}{2.0ex}$

\begin{lemm}\label{L3} If $V \cong {\cal A}(\a)$ or ${\cal B}(\a)$ as Vir-module,
then $V \cong {\cal A}(\a, 0)$ or ${\cal B}(\a, 0)$.
\end{lemm}

\noindent{\bf Proof.} If $V \cong {\cal A}(\a)$, then
$x_nv_i=(i+n)v_{n+i}$ if $t\ne 0$, $x_nv_0=n(n+a)v_{n}$ for some
$a\in\F$. We can deduce that $f(m, i)=0$ for all $n, i\in\z$.
Therefore $V\cong{\mathcal A}(a, 0)$.

If $V \cong {\cal A}(\a)$, then $x_nv_i=tv_{n+i}$ if $t\ne -n$,
$x_nv_{-n}=-n(n+a)v_{0}$,  for some $a\in{\mathbb F}$.

We can deduce that $f(m, i)=0$ for all $n, i\in\z$. Then $V$ is
isomorphic to ${\mathcal  B}(a, 0)$. \hfill
$\rule[-.23ex]{1.0ex}{2.0ex}$

\vskip10pt

From Lemma \ref{L2} and Lemma \ref{L3} we have
\begin{theo}\label{T3} Suppose  that  V is a nontrivial irreducible
weight ${\cal L}$-module with weight multiplicity one. Then we
have $V\simeq {\cal A}_{\a,\;\b,\;0}$ or $V\simeq {\cal
A}_{0,0,0}'$ for some $\a,\b\in\c$. Meanwhile, the three kinds
modules listed in the before of Lemma \ref{L1} are all
indecomposable weight ${\cal L}$-module with weight multiplicity
one.
\end{theo}

\section{Uniformly bounded irreducible weight modules
}

In this section,  we assume that $V$ is a uniformly bounded
nontrivial irreducible weight module over ${\cal L}$. So there
exists $\a\in\c$ such that $\hbox{Supp}(V)\subset \a+\z$. From
representation theory of $\Vir$, we have $C=0$ and $\dim
V_{\a+n}=p$ for all $\a+n \neq 0$. If $\a \in \z$, we also assume
 that $\a=0$.

Consider $V$ as a $\Vir$-module. We have a $\Vir$-submodule
filtration
    $$0=W^{(0)}\subset W^{(1)} \subset W^{(2)}\subset \cdots \subset W^{(p)}=V,$$
where $W^{(1)}, \cdots ,W^{(p)}$ are $\Vir$-submodules of $V$, and
the quotient modules \break $W^{(i)}/W^{(i-1)}$ have weight
multiplicity one for all nonzero weights.

Choose $v_n^1, \cdots, v_n^p \in V_{\a+n}$ such that the images of
 $v_n^i+W^{(i-1)}$ form a basis of $(W^{(i)}/W^{(i-1)})_{\a+n}$ for all $\a+n\neq 0$. We may suppose that
 $$x_i
 (v_n^1, \cdots ,v_n^p)=
 (x_iv_n^1, \cdots ,x_iv_n^p)=(v_{n+i}^1, \cdots ,v_{n+i}^p)A_{i,n},$$ where $A_{i,n}$
 are upper triangular $p\times p$ matrices, and $A_{i,n}(j,j)=\a+n+i\b_j $.
 Denote $$I(i)
 (v_n^1, \cdots ,v_n^p)=(v_{n+i}^1, \cdots ,v_{n+i}^p)F_{i,n}, \eqno (4.1)$$
where $F_{i,n}$
 are $p\times p$ matrices.

 The Lie brackets give
\vskip -.5cm
$$F_{i,j+n}F_{j,n}-F_{j,i+n}F_{i,n}=0,\eqno (4.2)$$
\vskip -.5cm
$$A_{i,j+n}A_{j,n}-A_{j,i+n}A_{i,n}=(j-i)A_{i+j,n},\eqno (4.3)$$
\vskip -.5cm
$$A_{i,j+n}F_{j,n}-F_{j,i+n}A_{i,n}=(j-i)F_{i+j,n}+{1\over 12}\delta_{i,-j}(i^3-i)C_1I_p,\eqno (4.4) $$
where the last three formulas have the restriction
$(\a+n)(\a+n+i)(\a+n+j)(\a+n+i+j)\ne0$. We shall denote the
$(i,j)$-entry of a matrix $A$ by $A(i,j)$.
\begin{lemm} \label{L41} If all nontrivial irreducible sub-quotient
$\Vir$-modules of V are isomorphic to ${\cal A}_{0,0}'$,
 then $V \simeq {\cal A}_{0,0,0}'$.
\end{lemm}
\noindent{\bf Proof.} Now we can suppose that $\dim (W^{(1)})_0\le
1$ (If $V$ contains a trivial submodule $\c v _0$, then the
span$\{u_k^1=I(k)v_0|k\in\z\}$ is a Vir-submodule, which can be
chosen as $W^{(1)}$).

 \noindent{\bf Claim}. {\it The $(k,1)$-entry $F_{j,n}(k,1)=0$ for all $k\geq 2$, $n \neq 0$ and
 $j+n \neq 0$.}

 \noindent{\it Proof of Claim.} Suppose  that we have $F_{j,n}(k,1)=0$ for all $k\geq k_0+1(k_0\geq 2)$, $n \neq 0$ and
 $j+n \neq 0$. We only need prove that $F_{j,n}(k_0,1)=0$ for all $n \neq 0$ and
 $j+n \neq 0$.

The $(k_0,1)$-entry of (4.4) gives
 $$(n+j)F_{j,n}(k_0,1)-F_{j,i+n}(k_0,1)n=(j-i)F_{i+j,n}(k_0,1),\,\,\,\hbox{if}\,\,\, n\ne 0,-i,-j,-i-j.\eqno(4.5)$$
Letting $j=1$ in (4.5), we have the $(k_0,1)$-entry
 $$(1-i)F_{i+1,n}(k_0,1)=(n+1)F_{1,n}(k_0,1)-nF_{1,i+n}(k_0,1),\,\,\,{\hbox{if}}\,\,\,n \neq 0,-1,-i,
 -1-i,$$
i.e.,
 $$(2-j)F_{j,n}({k_0},1)=(n+1)F_{1,n}({k_0},1)-nF_{1,n+j-1}({k_0},1),\,\,\,{\hbox{if}}\,\,\,n \neq 0,-1,-j+1,
 -j.\eqno(4.6)$$
 Letting $j=2$ in (4.6),  we have
$$nF_{1,n+1}({k_0},1)=(n+1)F_{1,n}({k_0},1),\,\,\,{\hbox{if}}\,\,\,n \neq
0,-1,-2.$$
 Hence $F_{1,n}({k_0},1)=nF_{1,1}({k_0},1)$
for all $n\geq 1$, and
$F_{1,n}({k_0},1)=-\frac{n}{2}F_{1,-2}({k_0},1)$ for all $n\le
-2$.

Moreover, by (4.5) we have
$$F_{m,n}({k_0},1)=nF_{1,1}({k_0},1), \ \ \  n\geq 1, m\ne -n-1, \eqno(4.a)$$
and
$$F_{m,n}({k_0},1)=-n/2F_{1,-2}({k_0},1), \ \ \  n\le -2, m\ne -n-1.\eqno(4.b)$$

Suppose  that $F_{1,1}({k_0},1)\neq 0$, then $F_{n,1}({k_0},1)\neq
0$ for all $n\ge 1$ by (4.a), then the span
$u_{n+1}=I(n)v_1^1+W^{(k_0-1)}\in W^{(k_0)}/W^{(k_0-1)}, n\in\z$,
is a nontrivial Vir-submodule. Moreover
$x_mu_{n+1}=x_mI(n)v_1^1+W^{(k_0-1)}=(n-m)I(m+n)v_1^1+W^{(k_0-1)}=(n-m)u_{m+n+1}$,
i.e. the Vir-module $V$ has a nontrivial submodule not isomorphic
to ${\cal A}_{0,0}'$, contradicting the assumption in the lemma.

Hence
$$F_{1,1}({k_0},1)=0.$$ Similarly we have
$$F_{1,-2}({k_0},1)=0.$$

Applying these to (4.6) we deduce that $F_{i,n}({k_0},1)=0$ for
all $n \neq 0, -1,-i,-i+1$. Letting $i=-n-1$ in (4.5) for suitable
$n$ we deduce that $F_{i,-1}({k_0},1)=0$, and letting $n=-i+1$ in
(4.5) for suitable $i$ we deduce that $F_{i,-i+1}({k_0},1)=0$.
 So we have proved this Claim.

This claim ensures that $I(i)v^1_j\in W^{(1)}$ if $j(i+j)\ne0$.
Consider the action of ${\cal L}$ on $W^{(1)}$. By the same
argument as in Lemma \ref{L3} we obtain that $F_{j,n}(1,1)=0$ for
all $j+n \neq 0$ and $n\neq 0.$

By computing the actions
$[x_i,I(j)]v^1_{-i-j}=(j-i)I(i+j)v^1_{-i-j}$ we deduce that \break
$\dim \sum_{j\in \z }\c I(j)v^1_{-j} \leq 1$. It is clear that
$I(k)I(j)v^1_{-j}=I(j)I(k)v^1_{-j}=0$ for all $k\ne j$, that
$(j-2k)I(j)I(j)v^1_{-j}=[x_k,I(j-k)]I(j)v^1_{-j}=0$ for many
suitable $k$. Hence all weight spaces of $W=U({{\cal L}})v_{1}^p$
is one dimensional. Combining with Theorem \ref{T3}, we have
proved this lemma. \hfill \qed

\begin{lemm} \label{L42}If any nontrivial irreducible sub-quotient
$\Vir$-module of V is  isomorphic to  ${\cal A}_{\a,0}$, where $\a
\notin \z$, then we have $V\simeq {\cal A}_{\a,0,0}$.
\end{lemm}
\noindent{\bf Proof.} We use the same notations and similar
discussions as in the proof of Lemma 4.1. Suppose  that we have
$F_{j,n}(k,1)=0$ for all $k\geq k_0+1(k_0\geq 2)$, $n\ne 0$ and
$j+n\ne 0$. We first want to prove that $F_{j,n}(k_0,1)=0$ for all
$n $ and
 $j$.

The $(k_0,1)$-entry of (4.4) gives
 $$(\a+n+j)F_{j,n}(k_0,1)-F_{j,i+n}(k_0,1)(\a+n)= (j-i) F_{i+j,n}(k_0,1).\eqno(4.7)$$
Letting $j=1$ in (4.7), we have the
 $$ (1-i)F_{i+1,n}(k_0,1)=(\a+n+1)F_{1,n}(k_0,1)-(\a+n)F_{1,i+n}(k_0,1),$$
i.e.,
 $$ (2-j)F_{j,n}({k_0},1)=(\a+n+1)F_{1,n}({k_0},1)-(\a+n)F_{1,n+j-1}({k_0},1).\eqno(4.8)$$
Letting $j=2$ in (4.8),  we have
 $0=(\a+n+1)F_{1,n}({k_0},1)-(\a+n)F_{1,n+1}({k_0},1)$.
Hence $$F_{1,n}({k_0},1)=\frac{\a+ n}{\a}F_{1,0}({k_0},1)
 ,\,\,\,\forall\,\,\,n\in\z.$$
Applying to (4.8) we obtain that
$$F_{j,n}({k_0},1)=\frac{\a+n}{\a}F_{1,0}(k_0,1),\,\,\,\forall\,\,\,j, n\in\z.\eqno(4.9)$$
Suppose  that $F_{1,0}(k_0,1)\ne0$. By re-scalaring
$\{v_i^{k_0}|i\in\z\}$ we may assume that $$
F_{1,0}(k_0,1)=\a.\eqno(4.10)$$

\noindent{\it Case 1: $k_0\geq 3$.}

{\it Case 1.1: $W^{(k_0)}/W^{(k_0-2)}$  is decomposable  over
$\Vir$.}
\medskip
In this case we can suitable choose $\{v^k_j|k,j\in\z\}$ so that
besides (4.9) we also have
$$F_{j,n}({k_0-1},1)=\frac{\a+n}{\a}F_{1,0}(k_0-1,1),\,\,\,\forall\,\,\,j, n\in\z.\eqno(4.11)$$
If $F_{1,0}(k_0 ,1)\ne0$, we know that $I(1)v^1_0 \;\mod
W^{(k_0-2)} ,v^{k_0}_1 \;\mod W^{(k_0-2)}$ are linearly
independent, and that $I(1)v^1_0 \;\mod W^{(k_0-2)} ,v^{k_0-1}_1
\;\mod W^{(k_0-2)}$ are linearly independent. Then we can
re-choose $W^{(k_0-1)}$ and $\{v^{k_0-1}_j| j\in\z\}$ such that
$v^{k_0-1}_1=I(1)v^1_0.$ Then $F_{1,0}({k_0},1)=0$, furthermore
 $F_{j,n}({k_0},1)=0 ,$ for all $j, n\in \z.$

{\it Case 1.2:    $W^{(k_0)}/W^{(k_0-2)}$ is indecomposable  over
$\Vir$.}  From Theorem 2.3, we need consider two subcases. \vskip
5pt

{\it Case 1.2.1:  $A_{i,n}(k_0-1,k_0)=-i$ for all $i, n \in \z$.}

 Using (4.9) and (4.10), from the $(k_0-1,1)$-entry of (4.4), we
 obtain
$$(\a+n+j)F_{j,n}(k_0-1,1)-i(\a+n)-F_{j,i+n}(k_0-1,1)(\a+n) =(j-i)F_{i+j,n}(k_0-1,1).\eqno(4.12)$$
 Letting $j=1$ and $i=-1$ in (4.12), we obtain that
$$(\a+n+1)F_{1,n}(k_0-1,1)=(\a+n)F_{1,n-1}(k_0-1,1)+2F_{0,n}(k_0-1,1)-(\a+n).\eqno(4.13)$$
Letting $i=j=1$ and $j=i=2$ in (4.12) respectively, we obtain that
  $$F_{1,1+n}(k_0-1,1)={\a+n+1\over \a+n}F_{1,n}(k_0-1,1)-1.\eqno(4.14)$$
$$F_{2,2+n}(k_0-1,1)={\a+n+2\over\a+n}F_{2,n}(k_0-1,1)-2.\eqno(4.15)$$

From (4.13) and (4.14) we have
$$F_{1,n}(k_0-1,1)=F_{0,n}(k_0-1,1)-1/2.\eqno(4.16)$$
$$(\a+n)F_{0,n+1}(k_0-1,1)=(\a+n+1)F_{0,n}(k_0-1,1)-\a-n-1/2.\eqno(4.17)$$
Letting $j=0$ and $i=2$ in (4.12) and using (4.17), we obtain that
$$F_{2,n}(k_0-1,1)=F_{0,n}(k_0-1,1)+\a+n-1/2.\eqno(4.18)$$
Combining (4.15), (4.17) and (4.18), we have $0=-2(\a+n)$, it is
contradiction. Then $F_{1,0}({k_0},1)=0$, furthermore
 $F_{j,n}({k_0},1)=0 \,\,\forall\,\,j, n\in \z.$

{\it Case 1.2.2:

$A_{\pm1,n}(k_0-1,k_0)=0,A_{\pm2,n}(k_0-1,k_0)=\pm
\frac{1}{(\a+n\pm 1)(\a+n\pm2)}.$}

Again using (4.9) and (4.10), from the $(k_0-1,1)$-entry of (4.4),
we obtain
$$(\a+n+j)F_{j,n}(k_0-1,1)+A_{i,j+n}(k_0-1,k_0)(\a+n)-F_{j,i+n}(k_0-1,1)(\a+n)$$
$$=(j-i)F_{i+j,n}(k_0-1,1).\eqno(4.19)$$ Letting $j=1$ and $i=-1$, we obtain that
$$(\a+n+1)F_{1,n}(k_0-1,1)+A_{-1,1+n}(k_0-1,k_0)(\a+n)-F_{1,n-1}(k_0-1,1)(\a+n),$$
$$=2F_{0,n}(k_0-1,1).\eqno(4.20)$$
i.e.
$$(\a+n+1)F_{1,n}(k_0-1,1)=F_{1,n-1}(k_0-1,1)(\a+n)+2F_{0,n}(k_0-1,1), \,\,\forall\,\,j, n\in \z.\eqno(4.21)$$
Letting $i=j=1$, we obtain that
$$(\a+n+1)F_{1,n}(k_0-1,1)-F_{1,1+n}(k_0-1,1)(\a+n)=0.\eqno(4.22)$$
Combining  (4.21) and (4.22), we deduce that
$$F_{1,n}(k_0-1,1)=F_{0,n}(k_0-1,1).\eqno(4.23)$$
Letting $i=j=2$, and $j=2,i=-1$ in (4.19) respectively, we obtain
that
$$(\a+n+2)F_{2,n}(k_0-1,1)+\frac{\a+n}{(\a+n+3)(\a+n+4)}-F_{2,2+n}(k_0-1,1)(\a+n)=0.\eqno(4.24)$$
$$(\a+n+2)F_{2,n}(k_0-1,1)-F_{2,n-1}(k_0-1,1)(\a+n)=3F_{1,n}(k_0-1,1).\eqno(4.25)$$
Combining  (4.24) and (4.25), we deduce that
$$F_{2,n+1}(k_0-1,1)={3\over \a+n}F_{1,n}(k_0-1,1)-{\a+n\over 3(\a+n+2)(\a+n+3)}.\eqno(4.26)$$
Combining  (4.24) and (4.26), we deduce that
$${6\over \a+n}F_{1,n}(k_0-1,1)={\a+n-1\over 3(\a+n+1)}-{(\a+n)(\a+n-2)\over 3(\a+n+3)(\a+n+4)}.\eqno(4.27)$$
It is contradict to (4.22). Then $F_{1,0}({k_0},1)=0$, furthermore
$F_{j,n}({k_0},1)=0 \,\,\forall\,\,j, n\in \z.$

{\it Case 2:  $k_0=2$.}

{\it Case 2.1:  $V$ is decomposable over $\Vir$. }

From the established Case 1 we may assume that $V=W^{(2)}$. Note
that $A_{i,n}(1,2)$ $=0$. The $(1,1)$-entry of (4.4) gives
$$(\a+j+n)F_{j,n}(1,1)-F_{j,i+n}(1,1)(\a+n)=(j-i)F_{i+j,n}(1,1)+{1\over 12}\delta_{i,-j}(i^3-i)C_1.\eqno(4.28)$$

As the calculation as in Lemma \ref{L1}, we obtain $C_1=0$, $F_{j,
n}(1,1)=d_{11}(\a+n)$ for some $d_{11}\in\c$ and for all $j,
n\in\z$.

Since $W^{(2)}$ is decomposable, by symmetry of (4.28) we have
$F_{j,n}(k,l)=d_{kl}(\a+n)$ for some $d_{kl}\in\c$, for all $j,
n\in\z$ and $k,l=1, 2$.

Thus $F_{j,n}=(\a+n)\left[\begin{array}{cc} d_{11} & d_{12}\\
d_{21} & d_{22}\end{array}\right]$.

By (4.2), we have $F_{j,n}=0$, so $V$ is an decomposable ${\cal
L}$-module and it is a contradiction.

{\it Case 2.2:  $W^{(2)}$ is indecomposable  over $\Vir$.}
\medskip The argument  is exactly the same as in Case 1.2.2. We do
not repeat it.

So far we have proved that $F_{j,n}(k,1)=0$ for all $n ,j\in\z$
and $k\ne1$. Thus $W^{(1)}$ is an ${\cal L}$-submodule which must
be $V$. Combining with Lemma \ref{L3}, we have proved this lemma.
 \hfill \qed

Denote by $(W^{(i)}/W^{(i-1)})'$ the unique nontrivial
sub-quotient $\Vir$-module of\break $W^{(i)}/W^{(i-1)}$. For any
$x,y \in \c$, define $x\nprec y$ if $y-x \notin \N $.

\begin{lemm}\label{L43} The module $V$ carries a filtration
$\{W^{(1)},W^{(2)}, \cdots ,W^{(p)}\}$ with \break
$(W^{(i)}/W^{(i-1)})'\simeq V'(\a,\b_i)$ as $\Vir$-modules, where
$ \b_i \nprec \b_j$ for all $i<j$.
\end{lemm}
\noindent{Proof}: We start with the filtration at the beginning of
this section. The statement is true if $p=1$, or $2$ (use Theorem
2.2 and $V'(\a,0)\simeq V'(\a,1)$). Now we consider $p>2$.

Suppose that we do not have $ \b_{i} \nprec \b_{i+1}$ for some
$i$, say $ \b_{p-1} \nprec \b_{p}$, i.e., $\b_{p}-\b_{p-1}\in \N$.
 Consider $V/V^{(p-2)}$ ($p=2$
for this module). Then we can have a submodule $ X^{(p-1)}\supset
V^{(p-2)}$ such that $(W^{(p)}/X^{(p-1)})'\simeq V'(\a,\b_{p-1})$
and $(X^{(p-1)}/W^{(p-2)})'\simeq V'(\a,\b_{p})$.

By repeating this procedure several times  if necessary, then we
obtain the filtration required. \hfill $ $\qed

Now we are ready to classify all irreducible uniformly bounded
weight  modules over ${\cal L}$.

\begin{theo}\label{T41}  If V is a nontrivial
irreducible uniformly bounded weight  module over ${\cal L}$, then
V is isomorphic to $V'(\a,\b; 0)$ for some $\a,\b\in \c$.
\end{theo}
\noindent{\bf Proof.} By Lemma 4.3, the module V carries a
filtration $\{W^{(1)},W^{(2)}, \cdots ,W^{(p)}\}$ with
$(W^{(i)}/W^{(i-1)})'\simeq V'(\a,\b_i)$ such that $$ \b_i \nprec
\b_j \,\,\, {\hbox{ for}}\,\,\,   {\hbox{ all}}\,\,\,  i<j. \eqno
(4.31)$$

From Lemma \ref{L41} and Lemma \ref{L42} we can suppose that
$\b_1\ne 0, 1$.

Suppose $F_{j,n}(k,1)=0$ for all $j,n\in\z$, $k>k_0$, where
$k_0>1$ is a fixed integer. We need to show that
$F_{j,n}(k_0,1)=0$ for all $j,n\in\z$.

\noindent{\bf Claim 1.} {\it $F_{1,n}(k_0,1)=0$  except for
finitely many $n\in\z$.}

{\it Case 1: $\a \notin \z$}.

 In this case all the restrictions for (4.2)-(4.4) disappear. Then the $(k_0,1)$-entry of  (4.4) gives
 $$(\a+n+j+i\b_{k_0})F_{j,n}(k_0,1)-F_{j,i+n}(k_0,1)(\a+n+i\b_1)=(j-i)F_{i+j,n}(k_0,1).\eqno(4.32)$$
 Letting $j=1$ we obtain
$$(1-i)F_{i+1,n}(k_0,1)=(\a+n+1+i\b_{k_0})F_{1,n}(k_0,1)-(\a+n+i\b_1)F_{1,i+n}(k_0,1).\eqno(4.33)$$
Taking $i=1$   we have
$$
(\a+\b_{k_0}+n+1)F_{1,n}({k_0},1)=(\a+\b_1+n)F_{1,n+1}({k_0},1).\eqno(4.34)$$
 Letting $i=-1$ in (4.33), we have
 $$2F_{0,n}({k_0},1)=(\a-\b_{k_0}+n+1)F_{1,n}({k_0},1)-(\a-\b_1+n)F_{1,n-1}({k_0},1).$$
 So
$$2(\a+n+\b_{k_0})F_{0,n}({k_0},1)=(\b_1^2-\b_{k_0}^2+2\a+2n+\b_{k_0}-\b_1)F_{1,n}(k_0,1). \eqno(4.35)$$
By using  (4.32) with $i=1$, $ j=0$, we deduce that
 $$(\a+n+\b_{k_0})F_{0,n}(k_0,1)-F_{0,1+n}(k_0,1)(\a+n+\b_1)=-F_{1,n}(k_0,1).\eqno(4.36)$$
Combining  (4.35), (4.36) and (4.34), we deduce that
$$F_{1,n}(k_0,1)=0\eqno(4.38)$$ except for finitely many
$n\in\z$

 {\it Case 2:  $\a=0$}.

Since $W^{(1)}\simeq A_{0,\b_1}$ and $\b_1\ne0$, $1$, then
$\dim(W^{(1)})_0=1$ and we can have $v_0^1$. In this case, the
restrictions for (4.4) become $(n+j)(n+i+j)\ne0.$ The restrictions
for (4.32)-(4.37) are $(n+j)(n+i+j)\ne0$, $(n+1)(n+i+1)\ne0$,
$n(n+1)\ne0$, $(n+1)(n+2)\ne0$, $(n+1)(n+2)(n+3)(n+4)\ne0$ and
$(n+1)(n+2)(n+3)(n+4)\ne0$, respectively. Thus we  have (4.38)
with exceptions $n=-1,-2,-3,-4$ and possibly one more exception
$n=n_0$ (which comes from the computation of getting (4.38)).
Claim 1 follows.

\noindent{\bf Claim 2.} {\it There exists some $i_0 \in \z, i_0
\neq -1,0$, such that $I(i_0+1)W^{(1)} \subseteq W^{(k_0-1)}$. }

For any $j\in\z\setminus\{0\}$, set
$$S_j=\{ n\in \z | I(j)v_{n}^1 \nsubseteq W^{(k_0-1)}\}.$$
By Claim 1 we know that $|S_1|< +\infty$. Choose $i_0$ satisfying

    (a) $i_0 \neq -1,0$;

    (b) $i_0> \max\{|x-y| | x, y \in S_1\}+1$;

    (c) $-\a+(-1+\b_1)i_0 \notin S_1, i_0\b_1+\b_1-\a-1 \notin S_1$ (because $\b_1 \neq
    0,1$).
 \vskip 5pt
 It is clear that
    $$(i_0k+S_1)\cap S_1=\emptyset,\,\,\,\hbox{for }\,\,\, k\ne0.$$
 From
    $$(\a+n-i_0+i_0\b)I(1)v_{n}^1= I(1)x_{i_0}v_{n-i_0}^1=(i_0-1)I(i_0+1)v_{n-i_0}^1+x_{i_0}I(1)v_{n-i_0}^1,$$
    we have
$$I(i_0+1)v_{n-i_0}^1 \in W^{(k_0-1)}\,\,\,\hbox{if}\,\,\,n \notin S_1\cup( i_0+S_1).\eqno(4.39)$$
Since $x_{-i_0}I(i_0+1)v_{n-i_0}^1 \in W^{(k_0-1)}$ for all $n
\notin S_1\cup(     i_0+S_1)$, i.e.,
      $$(2i_0+1)I(1)v_{n-i_0}^1+I(i_0+1)x_{-i_0}v_{n-i_0}^1$$
      $$=(2i_0+1)I(1)v_{n-i_0}^1+(\a+n-i_0-i_0\b)I(i_0+1)v_{n-2i_0}^1 \in W^{(k_0-1)},\,\,\,\forall \,\,\,n \notin S_1\cup(
      i_0+S_1),$$
we deduce that $$(\a+n-i_0-i_0\b)I(i_0+1)v_{n-2i_0}^1 \in
      W^{(k_0-1)},\,\,\,\forall \,\,\,n \notin S_1\cup(i_0+S_1).\eqno(4.40)$$
      For any $n \in 2i_0+S_1$, from (b) we have $n \notin S_1$ and
      $n \notin i_0+S_1,$  and from (c), we have $\a+n-i_0-i_0\b_1\neq 0$.
      Applying this to (4.40) we obtain that
      $$I(i_0+1)v_{\a+n}^1 \in
      W^{(k_0-1)},\,\,\,\forall \,\,\,n \in S_1.$$
      Together with (4.39), we have
   $$I(i_0+1)v_{n-i_0}^1 \in W^{(k_0-1)},\,\,\,\forall \,\,\,n\notin
   S_1.\eqno(4.41)$$
From $[x_{-i_0-1},I(i_0+1)]=2(i_0+1)I(0)$,
      we have $x_{-i_0-1}I(i_0+1)v_{n-i_0}^1\equiv
      I(i_0+1)x_{-i_0-1}v_{n-i_0}^1 \;\mod\; W^{(k_0-1)}.$
      Hence
     $$(\a+n-i_0-(i_0+1)\b_1)I(i_0+1)v_{n-2i_0-1}^1 \in
      W^{(k_0-1)},\,\,\,\forall \,\,\,n\notin S_1.\eqno(4.42)$$
     For $n-i_0-1\in S_1$ we know that $n\notin S_1$, and  by (c) we also have $(\a+n-i_0-(i_0+1)\b_1)\ne0$. Thus
     $$I(i_0+1)v_{n-i_0}^1 \in
      W^{(k_0-1)},\,\,\,\forall \,\,\,n\in S_1.$$
      Therefore $S_{i_0+1}=\emptyset.$
      So we have proved this Claim.

 Noting that $\{ x_i, I(i_0)|i\in\z\}$ generate  ${{\cal L}}$, we
 know that ${{\cal L}}W^{(1)}\subset W^{(k_0-1)}$. By induction on
 $k_0$ we see that $W^{(1)}$ is an ${\cal L}$-submodule. From
Theorem 3.4 we complete the proof of the theorem. \hfill \qed

\section{ Classification of irreducible weight modules over
${\cal L}$ with finite-dimensional weight spaces}
\medskip
\begin{theo}\label{T51} Let  V be an irreducible weight module over
${\cal L}$ with all weight spaces finite-dimensional. If V is not
uniformly bounded, then V is either a highest weight module or a
lowest weight module.
\end{theo}
\noindent{Proof.} Consider  $V$ as a $\Vir$-module. Let $W$ be the
smallest $\Vir$-submodule of $V$  such that $V/W$ is a trivial
$\Vir$-submodule. Then $W$ contains no trivial quotient module.

Let $W'$ be the maximal trivial $\Vir$-submodule of $W$. Then
$W/W'$ contains no trivial submodule.

Since $\dim W'+\dim V/W$ is finite, then $\Vir$-module $W/W'$ is
not uniformly bounded. Now $W/W'$ satisfies the conditions in
Theorem 2.4. By using Theorem 2.4 to $W/W'$, we have some
nontrivial upper bounded (or lower bounded) $\Vir $-submodule of
$W/W'$, say, $W''/W'$ is a nontrivial upper bounded $\Vir
$-submodule of $W/W'$ (i.e., the weight set of $W''/W'$ has an
upper bound).

Denote $W''=M^+$, We know that $M^+$ is not uniformly bounded.

For any $j\in\N$, define $M^+(j)=\{v\in M^+|I(i)
v=0\,\,\forall\,\, i\ge j\}$, and let
$$M=\cup_{i\in\N} M^+(j).$$
It is easy to check that $x_i M^+(j)\subset M^+(j+|i|)$, i.e., $M$
is an   $\Vir $-submodule of $M^{+}$. Suppose that
$\hbox{Supp}(V)\subset \a+\z$ for some $\a\in\c$.

\medskip
\noindent{\bf Claim. $M\ne0$.}  Fix $\lambda_0 \in \a+\z.$ Since $
M^+$ is not uniformly bounded, we have some $0\neq \lambda_1 \in
\hbox{Supp}(M^+)$ with $\lambda_1<\lambda_0$ and  $\dim
(M^+)_{\lambda_1}
>\dim V_{\lambda_0}$. Hence
$I(\lambda_0-\lambda_1):(M^+)_{\lambda_1}\rightarrow
V_{\lambda_0}$ is not injective. Say $v=v_{\lambda_1} \in
M^+\setminus\{0\}$ with $I(i_0)v=0$ where
$i_0=\lambda_0-\lambda_1>0$. Since $v\in M^+$, there exists
$j_0>0$ such that  $x_{j}v=0$ for $j\ge j_0$ Then we have
$I(i)v=0$ for all $i\ge i_0+j_0$. Thus $v\in M$. Claim follows.

\medskip
 Let $\Lambda$ be the maximal weight of $M$, and
$v_{\Lambda}$ is one of the corresponding weight vectors. By the
definition of $M$,  there exists a nonnegative integer $i_0$ such
that $I(i)v_{\Lambda}= 0$ for $i>i_ 0$, and
$I(i_0)v_{\Lambda}\ne0$ if $i_0>0$. If $i_0=0$, then $v_{\Lambda}$
is a highest weight vector of the ${{\cal L}}$-module $V$, and we
are done. So we assume that $i_0>0$. From
$x_jI(i_0)v_{\Lambda}=(i_0-j)I(i_0+j)v_{\Lambda}+I(i_0)x_jv_{\Lambda}=I(i_0)x_jv_{\Lambda}=0$
for all $j>0$, we know that $I(i_0)v_{\Lambda}\ne0 $ is a  highest
weight vector over $\Vir $. So $I(i_0)v_{\Lambda}\in M^+$. From
$I(i)I(i_0)v_{\Lambda}=I(i_0)I(i)v_{\Lambda}=0$ for all $i>i_0$,
we know that $I(i_0)v_{\Lambda}\in M$, contradicting the choice of
$\Lambda$ and $i_0>0$. So we have proved this Theorem. \hfill $
$\qed

 Combining Theorems \ref{T41} and Theorem \ref{T51} we obtain the main result of this paper:

\begin{theo} If V is a nontrivial irreducible weight
module over ${\cal L}$ with finite dimensional weight spaces, then
V is isomorphic to $V'(\a,\b; 0)$ for some $\a,\b\in \c$, or  a
highest or lowest weight module.
\end{theo}

In \cite{LGZ}, the authors proved the following theorem.
\begin{theo}\label{T52}\cite{LGZ} Let M be an irreducible
weight $\L$-module. Assume that there exists $\lambda\in\c$ such
that ${\rm dim\,}M_\lambda=\infty$. Then ${\rm
Supp}(M)=\lambda+\z$, and for every $k\in\z$, we have ${\rm
dim\,}M_{\lambda+k}=\infty$.
\end{theo}

With Theorem \ref{T52}, we classified all irreducible weight
modules of the $W$-algebra $W(2, 2)$.

\begin{theo} Let $M$ be an irreducible weight
$\cal L$-module. Assume that there exists $\lambda\in\c$ such that
$0<{\rm dim\,}M_\lambda<\infty$. Then $M$ is a Harish-Chandra
module. Consequently, $M$ is either an irreducible highest or
lowest weight module or an irreducible module from the
intermidiate series. \end{theo}

\def\refname{

\end{document}